\documentclass[12pt]{article}
\usepackage[dvips]{graphics}
\usepackage[pdftex]{graphicx}
\usepackage{graphicx}
\usepackage{color}
\usepackage{epsfig}
\usepackage{amssymb}
\usepackage{amsmath}
\usepackage{color}

\newtheorem{theorem}{Theorem}

\newtheorem{corollary}[theorem]{Corollary}

\newtheorem{statement}[theorem]{Statement}

\newcommand{\kom}[1]{}
\renewcommand{\kom}[1]{{\bf [#1]}}
\definecolor{blau}{rgb}{0.1,0.0,0.9}

\newcounter{komcounter}
\numberwithin{komcounter}{section}

\begin{document}

\title{{Analytical approximations of Lotka-Volterra integrals}}
\author{Niklas L.P. Lundstr\"om and Gunnar J. S\"oderbacka}


\date{}


\maketitle

\begin{abstract}
In this note we derive simple analytical bounds for solutions of $x - \ln x = y -\ln y$ and use them for estimating trajectories following Lotka-Volterra-type integrals. We show how our results imply estimates for the Lambert $W$ function as well as for trajectories of general predator-prey systems, including e.g. Rosenzweig-MacArthur equations. 

\end{abstract}

\section{Introduction}

\setcounter{theorem}{0}
\setcounter{equation}{0}

Consider the Lotka-Volterra system
\begin{align*}
&\frac{\mathrm dS}{\mathrm dt} = a S - b S X, \notag\\
&\frac{\mathrm dX}{\mathrm dt} = c S X - d X,
\end{align*}
where the non-negative variable $S = S(t)$ represents prey and the non-negative variable $X = X(t)$ predator biomass,
and $a$, $b$, $c$ and $d$ are positive constants.
By introducing the nondimensional quantities
%
$
x(\tau) = \frac{b}{a} X(t), \, s(\tau) = \frac{c}{d} S(t), \, \tau = a t, \, \alpha = \frac{d}{a}
$
%
we obtain the equivalent system
\begin{align}\label{eqmain}
&\frac{\mathrm ds}{\mathrm d\tau} = (1-x)s, \notag\\
&\frac{\mathrm dx}{\mathrm d\tau} = \alpha x (s-1),
\end{align}
which, by the corresponding phase-plane equation, has its trajectories on level curves of
\begin{align}\label{eq:Lyap-LT}
V(x,s) = \frac{1}{\alpha}\left( x - \ln x \right) + s - \ln s.
\end{align}
In this note we derive analytical bounds for solutions of the Lotka-Volterra integral \eqref{eq:Lyap-LT}.
Indeed, we prove that the solution $x < 1$ of the equation
\begin{align}\label{eq:maineq}
x - \ln x = y - \ln y, \quad \text{where} \quad y > 1,
\end{align}
satisfies the relation $x = z Y$, where $Y = y e^{-y}$ and $1 < z < e$ and $z = z(y)$ is a decreasing function of $y$.
We also prove that the inequalities $1 < z_1 < z < z_2 < z_0 < e$, where the $z_i$s are explicit functions of $Y$, hold for $z$, see Theorem \ref{le:julelemma} in Section \ref{sec:proof}.

\begin{figure}[ht!]
\begin{center}
\includegraphics[width=7.5cm,height=7cm]{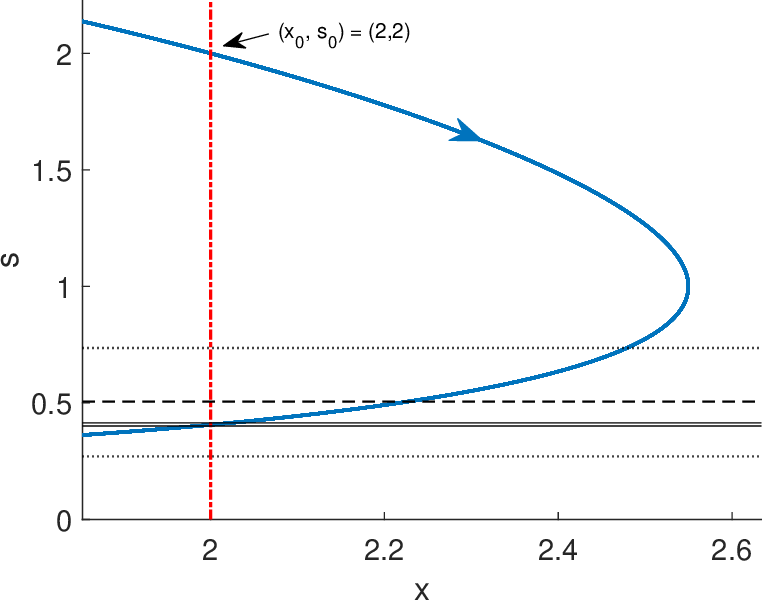}
\hspace{5mm}
\includegraphics[width=7.5cm,height=7cm]{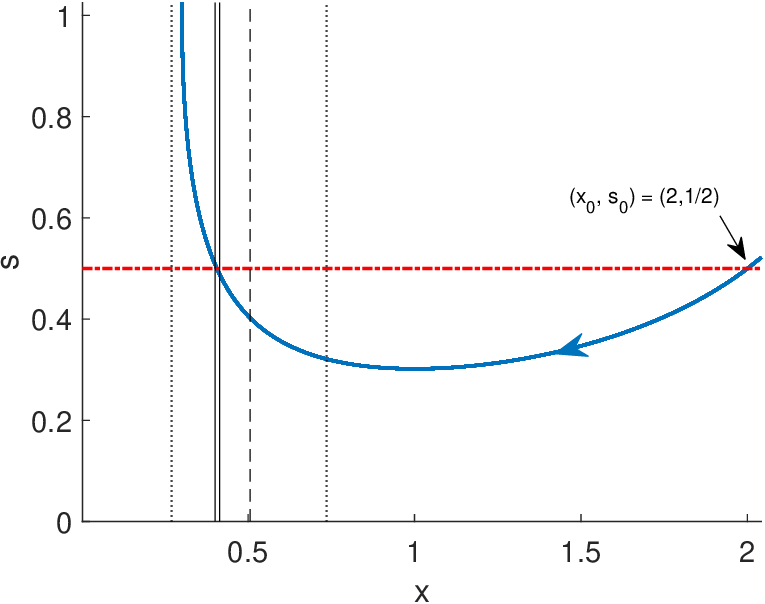}
\end{center}
\caption{Trajectory $T$ (solid blue curves) and the estimates of $T$ marked with straight lines in black: $z_1$ and $z_2$ (solid), $z_0$ (dashed) and $1$ and $e$ (dotted).
The red dashed-dot lines marks the sought after intersection level.
Estimates for $s$ in \eqref{eq:introny1} with $(x_0,s_0) = (2,2)$ (Left panel),
and for $x$ in \eqref{eq:introny2} with $(x_0,s_0) = (2,1/2)$ (Right panel). Here, $\alpha = 1$.}
\label{introfig}
\end{figure}

To apply the theorem, consider a trajectory $T$ of system \eqref{eqmain} with initial condition $(x_0, s_0)$ with $s_0 > 1$.
If we aim for an estimate of the $s$-value for the next intersection of $T$ with the line $x = x_0$ at, say $(x_0, s_1)$,
then we intend to estimate the solution of $V(x_0,s_0) = V(x_0,s_1)$,
where $V$ is given by \eqref{eq:Lyap-LT}, i.e. the solution $s_1 < 1 < s_0$ of
$$
s_1 - \ln s_1 =  s_0 - \ln s_0.
$$
According to the theorem we find
$$
s_1 = z(s_0) s_0 e^{-s_0}
$$
in which the function $z$ is decreasing and can be estimated so that
\begin{align}\label{eq:introny1}
1 < z_1 < \frac{s_1}{s_0} e^{s_0} < z_2 < z_0 < e
\end{align}
where $z_i = z_i(s_0)$, $i = 0,1,2$.
Likewise, if trajectory $T$ starts at $(x_0, s_0)$ with $x_0 > 1$,
then we may find an estimate of the next intersection of $T$ with the line $s = s_0$ through the equation $V(x_0,s_0) = V(x_1,s_0)$,
which boils down to $x_1 - \ln x_1 = x_0 - \ln x_0$.
As above we conclude, for the solution $x_1 < 1 < x_0$,
\begin{align}\label{eq:introny2}
1 < z_1 < \frac{x_1}{x_0} e^{x_0} < z_2 < z_0 < e
\end{align}
where $z_i = z_i(x_0)$, $i = 0,1,2$.
Observe that any level of predator or prey can be chosen,
giving estimates for the next intersection of the trajectory at the same predator or prey level.
For example, if we wish to estimate the minimal predator biomass, $x_{\text{min}}$,
on a trajectory having maximal predator biomass $x_{\text{max}}$,
then we use the fact that both maximal and minimal predator biomass are attained on the isocline $x' = 0$, i.e. at the prey biomass $s_0 = 1$.
From \eqref{eq:introny2} we obtain
\begin{align*}
1 < z_1 < \frac{x_{\text{min}}}{x_{\text{max}}} e^{x_{\text{max}}} < z_2 < z_0 < e
\end{align*}
where $z_i = z_i(x_{\text{max}})$, $i = 0,1,2$.
Clearly, we can obtain the similar estimates for $s_{\text{min}}$ as function of $s_{\text{max}}$ using \eqref{eq:introny1} and $x_0 = 1$.
Figure \ref{introfig} shows the estimates for $s_1$ in \eqref{eq:introny1} with $(x_0,s_0) = (2,2)$,
and the estimates for $x_1$ in \eqref{eq:introny2} with $(x_0,s_0) = (2,1/2)$.

In Theorem \ref{le:julelemma_high} we refine our arguments and derive more accurate bounds than those in Theorem \ref{le:julelemma} by introducing higher order Pad\'e approximants in the constructions.

A literature survey shows extensive interest in analysing Lotka-Volterra-type systems and their generalizations.
To mention a few, we refer the reader to \cite{MR87, CV98, GS07,IDM20} and the references therein.
Estimates valid for small prey biomass such as \eqref{eq:introny1}
may be of importance when studying the predators hunting strategies, e.g.,
if there is a threshold of the prey level at which the predator chose to switch from their central prey and instead starts to feed on other sources.
We give further motivations and demonstrations on how our theorems can be used to derive estimates of trajectories to more general dynamical systems in Section \ref{sec:appl}.

The solution of equation \eqref{eq:maineq} can be written
$$
x = -W(-ye^{-y})
$$
in which $W$ denotes the principal branch of the Lambert $W$ function.
Therefore, our estimates in Theorem \ref{le:julelemma} and Theorem \ref{le:julelemma_high} imply bounds and approximations of this function, see Corollary \ref{cor:Lambert} and Corollary \ref{cor:Lambert_high}.
Besides population dynamics, the Lambert $W$ function arises in many areas such as
chemical and mechanical engineering,
materials science,
statistical mechanics,
crystal growth,
economy,
viscous flows and flow through porous media, see e.g \cite{BPLPCS00,hej11,hej14,ACL24} and references therein.

In the next section, we present the theorem and its proof.



\section{Statement and proof of main results}

\setcounter{theorem}{0}
\setcounter{equation}{0}

\label{sec:proof}

In this section we prove the following analytical estimates:

\begin{theorem}\label{le:julelemma}
The solution $x<1$ to the equation $x- \ln x = y - \ln y$, where $y>1$, satisfies the relation $x=z\, Y$,
where $Y=y\, e^{-y}$, $1<z<e$, and $z = z(y)$ is decreasing in $y$.
Moreover,
%
%
the inequalities $1 < z_1 < z < z_2 < z_0 < e$ hold for $z$,
where, for $i = 1, 2$,
\begin{align*}
z_i &=\frac{1-d_i\, Y-\sqrt{(1-d_i\, Y)^2 - 4 c_i\, Y}}{2c_i\, Y},
\quad  z_0 =\frac{1}{1-(e-1)Y},\notag\\
\quad d_i &= e-1-c_i\, e,
\quad c_1=\frac{e-2}{e-1},
\quad \text{and} \quad c_2=\frac{1}{e}.
\end{align*}
\end{theorem}

\begin{figure}[ht!]
\begin{center}
\includegraphics[width=7.5cm,height=6.5cm]{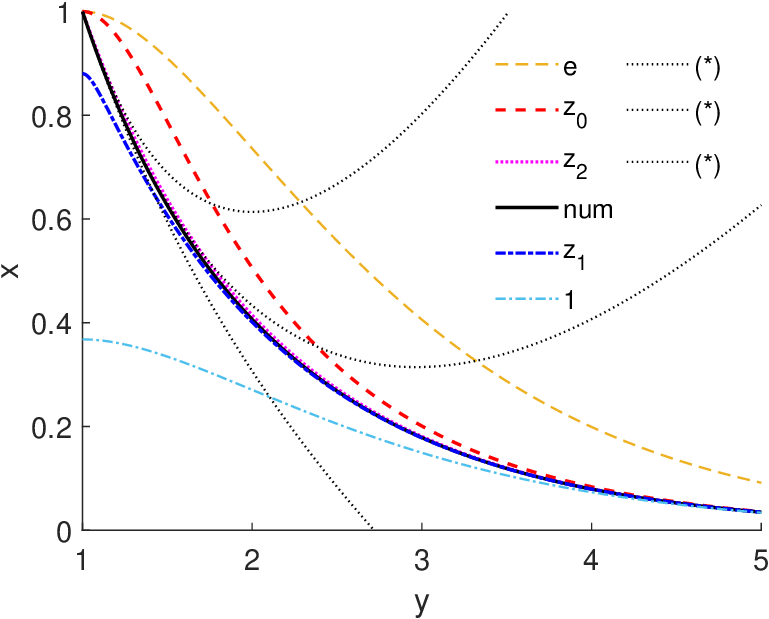}
\hspace{5mm}
\includegraphics[width=7.5cm,height=6.5cm]{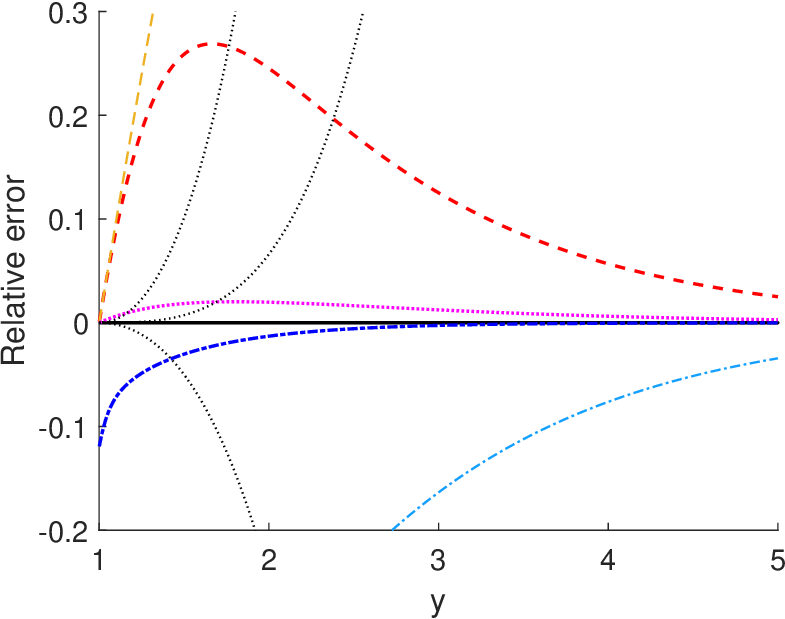}
.\vspace{5mm}
\includegraphics[width=7.5cm,height=6.5cm]{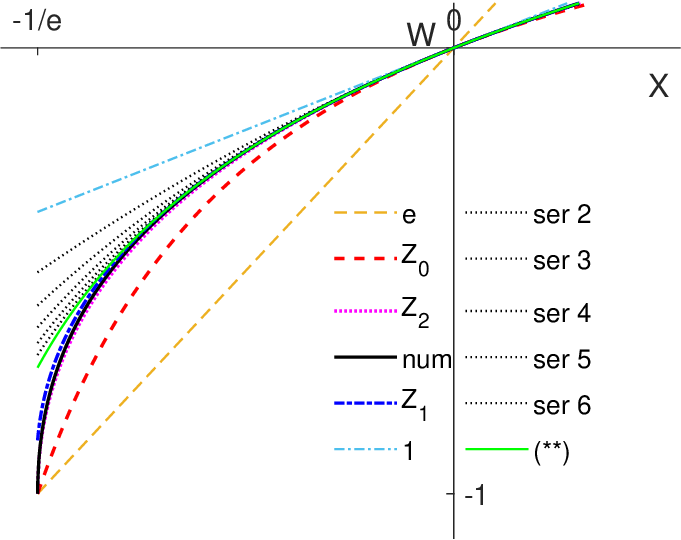}
\hspace{5mm}
\includegraphics[width=7.5cm,height=6.5cm]{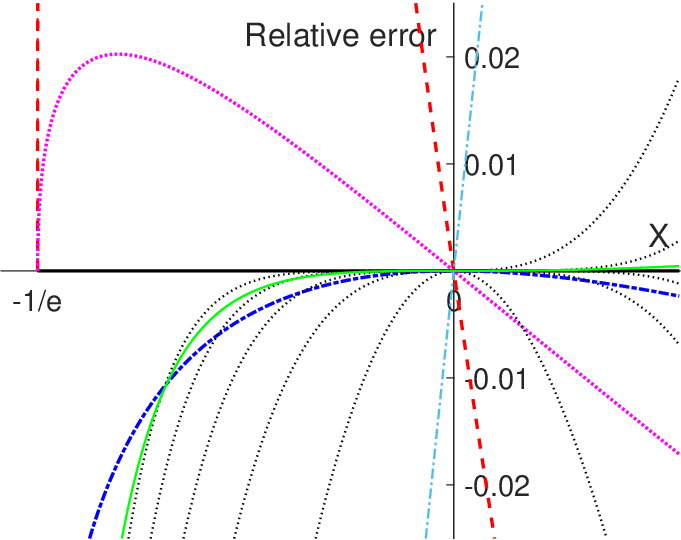}
\end{center}
\caption{
A numerical solution of the equation $x -\ln x = y - \ln y$ together with the five bounds in Theorem \ref{le:julelemma}, as well as the bounds in display \eqref{eq:S1} (upper left).
The Lambert $W$-function together with the bounds in Corollary \ref{cor:Lambert}, the upper bound in display \eqref{eq:wikibound} with $\Bar{y} = X + 1$ and the series approximation in display \eqref{eq:series} with 2, 3, 4, 5 and 6 terms (lower left).
Relative error (right panels).}
\label{xappro}
\end{figure}

\noindent
A numerical solution of the equation $x -\ln x = y - \ln y$ together with the five bounds in Theorem \ref{le:julelemma}, as well as the bounds in \eqref{eq:S1},
are plotted in Figure \ref{xappro} (upper left panel).
The upper right panel shows the relative error.

While the estimates in Theorem \ref{le:julelemma} are not impressively accurate, we emphasize their simplicity and the fact that in biological systems, the Lotka-Volterra integral constitutes already an approximation of a real system, motivating us to strive for simple expressions rather than higher precision.
Before giving the proof we also remark that any equation of type $x - a\, \ln x = y - a\, \ln y\, 0<x<a<y$
can be transformed into the equation of Theorem \ref{le:julelemma} by scaling
$x$ and $y$.

\bigskip

\noindent
{\it Proof of Theorem \ref{le:julelemma}}.
We begin by proving the first statement of the Theorem by
substituting $x=z\, Y$ into the equation of the theorem to obtain
$$
z Y - \ln z - \ln Y = y - \ln y
$$
and thus, since $Y = y e^{-y}$ we have,
\begin{equation}\label{z}
Z(z)=Y\quad \text{in which} \quad Z(z)=\frac{\ln z}{z}.
\end{equation}
This equation has a unique solution for $z$, where  $1<z<e$,
because $0<Y<1/e$ for $y>1$ and $Z$ is increasing in $z$ and $Z(1)=0$ and $Z(e)=1/e$.
Differentiation of formula (\ref{z}) with respect to $y$,
$$
\frac{\mathrm d}{\mathrm dy}\left( y e^{-y}\right) = \frac{\mathrm d}{\mathrm dy}\left( \frac{\ln z(y)}{z(y)}\right),
$$
gives
$$
z' = z^2\, e^{-y}\, \frac{1-y}{1-\ln z} < 0.
$$
Hence, $z$ is decreasing in $y$.

In order to get the estimates for $z$ we intend to replace $\ln z$ in (\ref{z}) by Pad\'e approximations built on the rational functions
\begin{align}\label{eq:pade_orig}
f_i(z)=\frac{z-1}{c_i \, z +d_i}, \quad i = 0, 1, 2, \quad  1 \leq z \leq e,
\end{align}
%
and then solve the remaining formulas for $z$.
Immediately $f_i(1) = 0$, and by demanding $f_i(e) = 1$ we obtain $d_i = e - 1 - c_i e$ for $i = 0, 1, 2$.
Taking $c_0 = 0$ makes $f_0$ a linear approximation and equating the derivatives of $\ln z$ and $f_2(z)$ at $e$ gives 
$c_2 = 1/e$.
Similarly, equating the derivatives of $\ln z$ and $f_1(z)$ at $1$ 
gives $c_1 = (e-2)/(e-1)$,
and we prove below that these choices imply the central inequalities
\begin{align}\label{eq:ineq:lotka}
f_0(z)<f_2(z)< \ln z < f_1(z) \quad \text{for} \quad  1 < z < e.
\end{align}

Denote by $z^*$ the solution to equation (\ref{z}).
Then $Z(z)<Y$ for $z<z^*$ and $Z(z)>Y$ for $z>z^*$.
Thus if $Y=\frac{f_i(z_i)}{z_i}<Z(z_i)$ then $z^*<z_i$ and if $Y=\frac{f_i(z_i)}{z_i}>Z(z_i)$ then $z^*>z_i$.
In order to prove \eqref{eq:ineq:lotka} we now consider the functions $g_i$ defined by
$$
g_i(z)=\ln z -f_i(z).
$$
Calculating the derivative of $g_i$ we get
\begin{equation*}
g_i'(z)=\frac{h_i(z)}{z\, (c_i\, z+d_i)^2}\quad\text{in which} \quad h_i(z)=c_i^2\, z^2 + (2d_i\, c_i -c_i -d_i)\, z + d_i^2.
\end{equation*}
We notice that $h_1$ is negative between 1 and $\frac{1}{(e-2)^2}$, positive between
$\frac{1}{(e-2)^2}$ and $e$,
and because $g_1(1)=g_1(e)=0$, we conclude
that $g_1(z)<0$ between 1 and $e$.
Thus $\frac{f_1(z)}{z}>Z(z)$ and because $z_1$ is the solution to
$\frac{f_1(z)}{z}=Y$ we get $z^*>z_1$.

Further, $h_2$ is positive between 1 and $e\, (e-2)^2$, negative between $e\, (e-2)^2$ and $e$ and because $g_2(1)=g_2(e)=0$,
we conclude that $g_2(z)>0$ between 1 and $e$.
Thus $\frac{f_2(z)}{z}<Z(z)$ and because $z_2$ is the solution to
$\frac{f_2(z)}{z}=Y$ we get $z^*<z_2$.

Finally, $h_0$ is positive between 1 and $e-1$, negative between $e-1$ and $e$ and because $g_0(1)=g_0(e)=0$,
we conclude that $g_0(z)>0$ between 1 and $e$.
Thus $\frac{f_0(z)}{z}<Z(z)$ and because $z_0$ is the solution to
$\frac{f_0(z)}{z}=Y$   we get $z^*<z_0$. $\hfill\Box$\\

Using higher order Pad\'e approximations in place of \eqref{eq:pade_orig} we next build the following tighter bounds:

\begin{theorem}\label{le:julelemma_high}
The inequalities $\tilde z_1 < z < \tilde z_i$, $i = 2,3$, hold for $z$ in Theorem \ref{le:julelemma},
where
\begin{align*}
\tilde z_i &=\frac{1 - 2a_i - d_i\, Y -\sqrt{(1-d_i\, Y)^2 - 4 Y(c_i - a_i(d_i + c_i))}}{2(c_i\, Y - a_i)},
\end{align*}
in which $d_i = e-1-c_i\, e + a_i (e-1)^2$, $i = 1, 2, 3$ and where 
\begin{align*}\label{koeffs_high}
a_1 = 1 - \frac{e}{(e-1)^2}, \qquad a_2 &= \frac{3 - e}{2(e - 1)(e - 2)}, \qquad  
a_3  = \frac{c_3\, e - 1}{e^2 - 1},
\notag\\
c_1 = e-1 -\frac{2}{e-1}, \qquad c_2 &= \frac{e^2 - 4e + 5}{2(e - 1)(e - 2)}, \qquad
c_3 = \frac{2e - (e-1)^2}{2 + (e-1)^2}.
\end{align*}
\end{theorem}

\noindent
{\it Proof of Theorem \ref{le:julelemma_high}}.
The argument is very similar to the second part of the proof of Theorem \ref{le:julelemma} instead of \eqref{eq:pade_orig} we estimate $\ln z$ with the higher order Pad\'e
$$
f_i(z)=\frac{z-1 + a_i(z-1)^2}{c_i \, z + d_i}, \quad i = 1, 2, 3,
\quad  1 \leq z \leq e.
$$
Solving the remaining expression of \eqref{z}, which is only a second order equation, gives the desired expression for $z$,
and equating $\ln$ with $f_i$ at $e$ immediately gives
$
d_i = e - 1 - c_i e + a_i (e - 1)^2.
$
%
%
We will show
$$
f_1(z) < \ln z < f_j(z), \quad j = 2,3, \quad 1 < z < e,
$$
by observing that the derivative of
$
g_i(z)=\ln z -f_i(z)
$
yields
\begin{equation*}
g_i'(z)=\frac{h_i(z)}{z\, (c_i\, z+d_i)^2}, \quad i = 1, 2, 3,
\end{equation*}
in which
$$
h_i(z) = -a_i c_i\, z^3 + (c_i^2 - 2a_i d_i)\, z^2  + (a_i c_i + 2 a_i d_i - d_i - c_i + 2 d_i c_i)\, z + d_i^2.
$$

Equating first derivatives of $\ln$ and $f_1$ at endpoints 1 and $e$ gives $a_1$ and $c_1$ building the lower bound $\tilde z_1$.
We notice that $h_1$ has three real roots, 1, $z_r \approx 1.66$ and $e$, is negative between 1 and $z_r$, positive between $z_r$ and $e$,
and because $g_1(1) = g_1(e) = 0$ we conclude that $g_1(z) < 0$ between 1 and $e$.
Thus $\frac{f_1(z)}{z}>Z(z)$ and because $\tilde z_1$ is the solution of $\frac{f_1(z)}{z}=Y$ we get $z^*>\tilde z_1$.

Equating first and second derivatives of $\ln$ and $f_1$ at $z = 1$ gives $a_2$ and $c_2$ building the upper bound $\tilde z_2$.
We notice that $h_2$ has three real roots, 1, 1, $z_r \approx 2.12$, is positive between 1 and $z_r$, negative between $z_r$ and $e$ and because $g_2(1)=g_2(e)=0$,
we conclude that $g_2(z)>0$ between 1 and $e$.
Thus $\frac{f_2(z)}{z}<Z(z)$ and because $\tilde z_2$ is the solution of $\frac{f_2(z)}{z}=Y$ we get $z^*<\tilde z_2$.

Equating first and second derivatives of $\ln$ and $f_1$ at $z = e$ gives $a_3$ and $c_3$ building the upper bound $\tilde z_3$.
We notice that $h_3$ has three real roots, $z_r\approx 1.296,e,e$, is positive between 1 and $z_r$, negative between $z_r$ and $e$ and because $g_3(1)=g_3(e)=0$
we conclude that $g_3(z)>0$ between 1 and $e$.
Thus $\frac{f_3(z)}{z}<Z(z)$ and because $z_3$ is the solution of
$\frac{f_3(z)}{z}=Y$ we get $z^*<\tilde z_3$. $\hfill\Box$\\

\begin{figure}[ht!]
\begin{center}
\includegraphics[width=7.5cm,height=6.5cm]{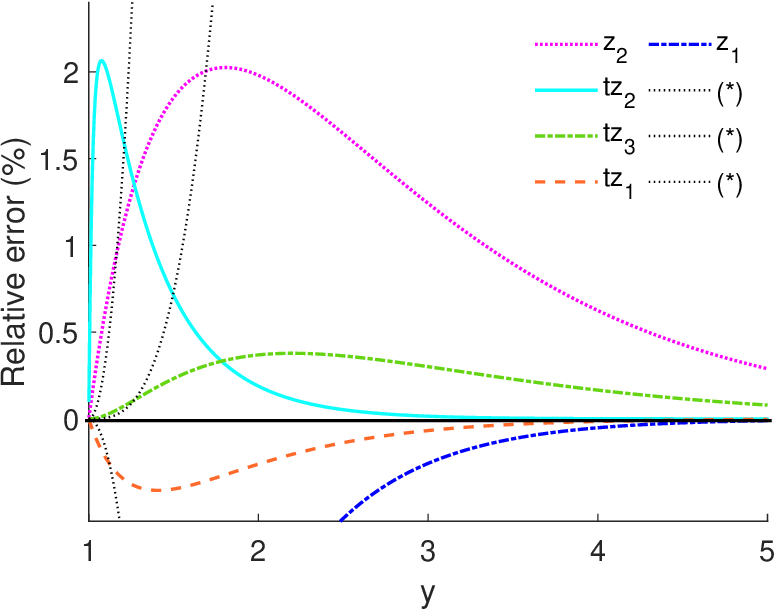}
\hspace{5mm}
\includegraphics[width=7.5cm,height=6.5cm]{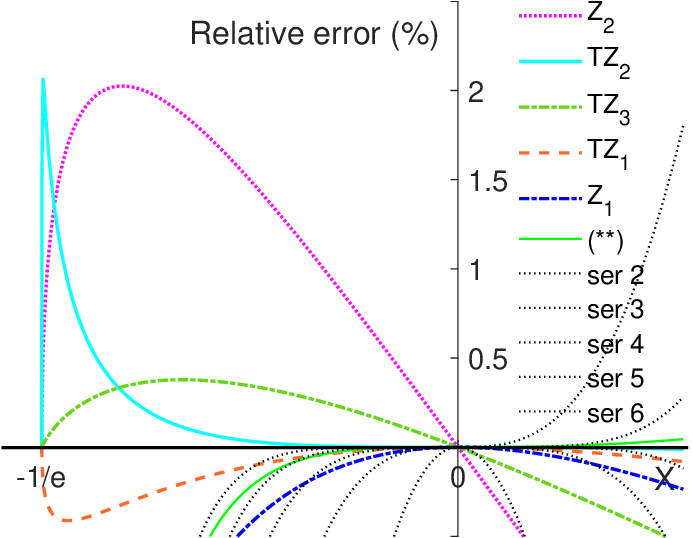}
\end{center}
\caption{Relative error of the bounds in Theorem \ref{le:julelemma_high}, the sharpest bounds from Theorem \ref{le:julelemma} and those given in display \eqref{eq:S1} (left panel).
Relative error of the bounds on the Lambert $W$-function in Corollary \ref{cor:Lambert_high},
the sharpest bounds from Corollary \ref{cor:Lambert}, the bound in display \eqref{eq:wikibound} with $\Bar{y} = X + 1$, and the series approximation in display \eqref{eq:series} with 2, 3, 4, 5 and 6 terms (right panel).
In the legend, $tz_i = \tilde z_i$ and $TZ_i = \widetilde{\mathcal{Z}_i}$.
}
\label{fig:high}
\end{figure}

Of the two upper bounds, they match better near the side where the derivatives are equated, naturally.
Figure \ref{fig:high} (left panel) shows the relative error of the bounds in Theorem \ref{le:julelemma_high}, the sharpest bounds from Theorem \ref{le:julelemma} and those given in \eqref{eq:S1}.

\subsection*{Implications for the Lambert $W$-function}

\begin{figure}[ht!]
\begin{center}
\includegraphics[width=7.5cm,height=7cm]{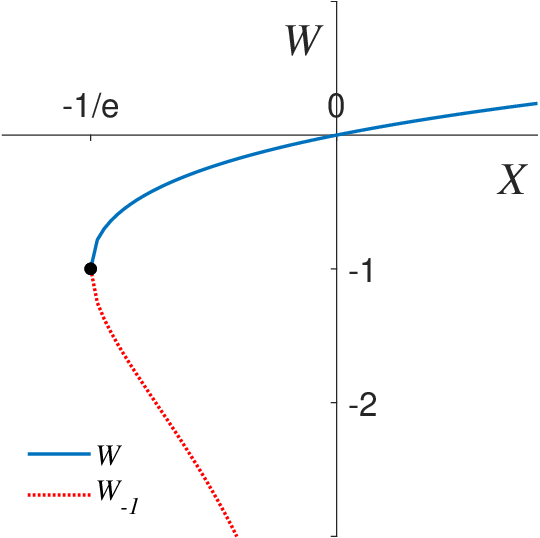}
\end{center}
\caption{The Lambert W function.}
\label{fig:Lambert}
\end{figure}

For real numbers $X$ and $u$ the equation
$$
u e^{u} = X
$$
can be solved for $u$ only if $X \geq -1/e$; gets $u = W(X)$ if $X \geq 0$ and the two values $u = W(X)$ and $u = W_{-1}(X)$ if $-1/e \leq X < 0$.
Here $W$ is the upper (principal) branch and
$W_{-1}$ the lower branch of the Lambert W function,
see Figure \ref{fig:Lambert}.

The equation $x-\ln{x} = y - \ln{y}$ , $x \in (0,1)$, $y\in(1,\infty)$ can be written
$$
-x e^{-x} = -y e^{-y} = -Y
$$
and thus $x = -W(-ye^{-y}) = -W(-Y)$.
On the other hand, following the notation in Theorem \ref{le:julelemma} we also have
$$
x = z\, ye^{-y} = -W(-ye^{-y}) = -W(-Y),
$$
and hence Theorem \ref{le:julelemma} gives estimates of the function $-W(-ye^{-y})$.
In our case $-Y \in (-e^{-1},0)$, $x \in (0,1)$ and hence we are in the upper (principal) branch.
We remark that the function $W(-ye^{-y})$ appears also in the classic problem of a projectile moving through a linearly resisting medium \cite{WW04, PY04, M05, S06}, and that several bounds for $W(-ye^{-y})$ were derived in \cite{S09}. For example \cite[Theorem 3.5 and Theorem 3.7]{S09} imply
\begin{align}\label{eq:S1}
2 \ln{y} - y < \sqrt{8(y-1 - \ln{y})} - y < W(-ye^{-y}) < \ln{y} - 1 \tag{$\star$}
\end{align}
whenever $y > 1$.
(In \cite{S09} the r.h.s. is $2\ln{y} - 1$, but their proof holds for \eqref{eq:S1} as well.)
Figure \ref{xappro} (upper panels) shows the function $x = - W(-ye^{-y})$ together with our estimates and the estimates in \eqref{eq:S1}.

%
%

Next, let $X = -Y$ and observe that
$$
W(X) = W(-Y) = -z Y = z X.
$$
Noticing also that $z = z(y)$ depends only on $Y = -X$, the estimates in Theorem \ref{le:julelemma} implies the following approximations of the Lambert $W$-function:

\begin{corollary}\label{cor:Lambert}
Let $W(X)$ be the principal branch of the Lambert $W$-function. 
Then
$$
\mathcal Z_1( X) \leq W( X) \leq \mathcal Z_2( X) \leq \mathcal Z_0( X)
$$
where, for $i=1,2$,
\begin{align*}
\mathcal Z_i(X) &=\frac{-1-d_i\, X+\sqrt{(1+d_i\, X)^2 + 4 c_i\, X}}{2c_i},
\quad \mathcal Z_0(X) = \frac{X}{1+(e-1)X},\\
\quad d_i &= e-1-c_i\, e,
\quad c_1=\frac{e-2}{e-1},
\quad \text{and} \quad c_2=\frac{1}{e}.
\end{align*}
\end{corollary}

\noindent
While the estimates in Corollary \ref{cor:Lambert} are not impressively accurate, we emphasize their simplicity
and that we will derive more accurate bounds for the Lambert $W$-function in Corollary \ref{cor:Lambert_high} below.
We proceed by comparing the bounds in Corollary \ref{cor:Lambert} with other simple estimates, e.g. the following upper bound from \cite{HH08}:
\begin{align}\label{eq:wikibound}
W(X) \leq \left(\frac{X + \Bar{y}}{1 + \ln\Bar{y}}\right)
\tag{$\star\star$}
\end{align}
valid for $X \geq - 1 / e$ where $\Bar{y} > 1 / e$ is a degree of freedom.
Moreover, the Taylor series of $W$ around 0 yields
\begin{align}\label{eq:series}
W(X) = \sum_{n=1}^\infty \frac{(-n)^{n-1}}{n!} X^n = X - X^2 + \frac{3}{2} X^3 - \frac{8}{3} X^4 + \frac{125}{24} X^5 - \dots, \tag{ser}
\end{align}
and an approximation with relative error less than $0.013\%$ can be found in \cite[eqs (7), (8) and (9)]{BPLPCS00}, to which we also refer the reader for an extensive list of applications for the Lambert W function.
Figure \ref{xappro} (lower panels) shows the Lambert $W$-function together with our estimates in Corollary \ref{cor:Lambert}, the upper bound \eqref{eq:wikibound} with $\Bar{y} = X + 1$ and the series approximation \eqref{eq:series} with 2, 3, 4, 5 and 6 terms (left), and the relative error (right).

In the same way as Theorem \ref{le:julelemma} yield Corollary \ref{cor:Lambert}, the tighter bounds in Theorem \ref{le:julelemma_high} implies the following tighter bounds of the Lambert $W$-function:

\begin{corollary}\label{cor:Lambert_high}
Let $W(X)$ be the principal branch of the Lambert $W$-function. 
Then
$$
\widetilde{\mathcal Z_1}( X) \leq W( X) \leq  \widetilde{\mathcal Z_i}(X), \quad i = 2,3,
$$
where
\begin{align*}
\widetilde{\mathcal Z_i}(X) &=\frac{2a_i\,-1-d_i\, X+\sqrt{(1+d_i\, X)^2 + 4 X(c_i - a_i(d_i + c_i))}}{2(c_i + a_i X^{-1})}
, \quad i=1,2,3,
\end{align*}
%
%
%
and where coefficients $a_i, c_i$ and $d_i$ are as in Theorem \ref{le:julelemma_high}.
\end{corollary}

Figure \ref{fig:high} (right panel) shows the relative error of the bounds on the Lambert $W$-function in Corollary \ref{cor:Lambert_high},
the sharpest bounds from Corollary \ref{cor:Lambert},  the bound in \eqref{eq:wikibound} with $\Bar{y} = X + 1$ and the series approximation \eqref{eq:series} with 2, 3, 4, 5 and 6 terms.


\section{Applications to more general predator-prey systems}

\setcounter{theorem}{0}
\setcounter{equation}{0}

\label{sec:appl}

%
%

Consider a general predator-prey system on the form
\begin{align}\label{eq:general}
\frac{\mathrm dS}{\mathrm dt} &= H(S) - q \varphi(S) X, \notag\\
\frac{\mathrm dX}{\mathrm dt} &= p \varphi(S) X - d X,
\end{align}
where the non-negative variable $S = S(t)$ represents the prey biomass,
the non-negative variable $X = X(t)$ represents the predator biomass,
$\varphi$ is non-decreasing,
$\varphi(0) = H(0) = 0$,
and parameters $p, q, d$ are positive.
Systems of type \eqref{eq:general} have been extensively studied through the last centuries, see e.g. \cite{C81,L89,L94,HS09,LS18,LS22} and the references therein.

Often, functions $H$ and $\varphi$ are defined by
\begin{align}\label{eq:RM}
H(S) = r S\left(1 -\frac{S}{K} \right) \quad \text{and}\quad \varphi(S) = \frac{S^n}{S^n + A},
\end{align}
where most commonly $n = 1$ or $n = 2$.
In case of \eqref{eq:RM} system \eqref{eq:general} is usually referred to as a Rosenzweig-MacArthur predator prey system.
If
\begin{align*}
H(S) = r S \quad \text{and}\quad \varphi(S) = S,
\end{align*}
then system \eqref{eq:general} boils down to the Lotka-Volterra equations \eqref{eqmain}.

We now intend to analyze the general system \eqref{eq:general} using our estimates in Theorem \ref{le:julelemma}.
Assume, without loss of generality, that $q = 1$.
The phase-plane equation yields
\begin{align*}
\frac{\mathrm dS}{\mathrm dX} = \frac{F(S) -  X}{X} \cdot\frac{\varphi(S)}{p \varphi(S) - d}
\end{align*}
where $F(S) = H(S) / \varphi(S)$.
Let us replace $F(S)$ by a constant $\overline F$ for the moment, and observe that then integrating gives
\begin{align*}
\int\left(p - \frac{d}{\varphi(S)}\right)dS  =  \int\left(\frac{{\overline F}}{X} - 1\right) dX
\end{align*}
and thus the system can, under reasonable assumptions on $\varphi$ and $H$, be analyzed by the generalized Lotka-Volterra integral
\begin{align*}
V_{\overline F}(X,S) = pS - d\int\frac{dS}{\varphi(S)}  +  X - {\overline F} \ln{X}.
\end{align*}
Moreover
\begin{align}\label{eq:japp1}
\nabla V_{\overline F} = \left(1 - \frac{\overline F}{X} , p - \frac{d}{\varphi(S)}\right)
\end{align}
and
\begin{align}\label{eq:japp2}
\frac{\partial V}{\partial t}
&= \left(p - \frac{d}{\varphi(S)}\right) S' + \left(1 - \frac{\overline F}{X} \right) X'
= \left(p\varphi(S) - d\right) \left(F(S) - {\overline F} \right).
\end{align}

\begin{figure}[ht!]
\begin{center}
\includegraphics[width=8cm,height=3.5cm]{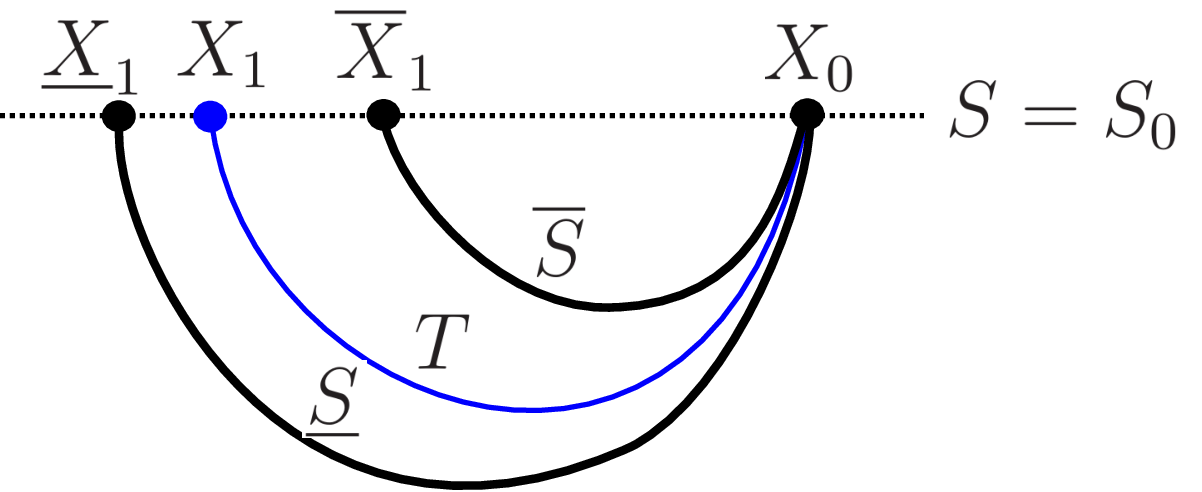}
\end{center}
\caption{Geometry in the construction of estimates.}
\label{fig:trapped}
\end{figure}

To proceed we consider a trajectory $T$ of system \eqref{eq:general} with initial condition $(X_0,S_0)$ where $X_0 > F(S)$ and $S_0$ satisfies $\varphi(S_0) < d/p$.
Since $X = F(S)$ gives an isocline with $S' = 0$ and
$\varphi(S) = d/p$ gives an isocline with $X' = 0$,
this means that both $S$ and $X$ are decreasing initially.
Suppose that, until $T$ intersects $S = S_0$ the next time, $T$ stays in a part of the state space where there exist positive $\underline{F}$ and $\overline{F}$ such that
\begin{align}\label{eq:ass_general}
\underline{F} < F(S) < \overline{F} \quad \text{and} \quad \varphi(S) < d/p.
\end{align}
It then follows from \eqref{eq:japp1} and \eqref{eq:japp2} that trajectory $T$ starting at $(X_0,S_0)$ remains trapped between the curves $\underline{S}$ and  $\overline{S}$, defined through
$$
V_{\underline{F}}(X_0,S_0) = V_{\underline{F}}(X,\underline{S}(X)) \quad \text{and} \quad
V_{\overline{F}}(X_0,S_0) = V_{\overline{F}}(X,\overline{S}(X)),
$$
see Figure \ref{fig:trapped}.
\noindent
Moreover, the `barriers' $\underline{S}$ and $\overline{S}$ are convex with min at
$X = \underline{F}$ and $X = \overline{F}$,
and intersect $S = S_0$ a second time at
$\underline{X}_1$ and $\overline{X}_1$, respectively.
For the next intersection of $T$ with $S = S_0$ at $(X_1,S_0)$ it necessarily holds that
$$
\underline{X}_1 < X_1 < \overline{X}_1,
$$
where
\begin{align*}
V_{\underline{F}}(X_0,S_0) = V_{\underline{F}}(\underline{X}_1,S_0) \quad \text{giving}\quad  X_0 - \underline{F} \ln{X_0} =  \underline{X}_1 - \underline{F} \ln{\underline{X}_1},
\end{align*}
and
\begin{align*}
V_{\overline{F}}(X_0,S_0) = V_{\overline{F}}(\overline X_1,S_0) \quad \text{giving}\quad  X_0 - \overline{F} \ln{X_0} = \overline{X}_1 - \overline{F} \ln{\overline{X}_1}.
\end{align*}
By setting $x = \frac{\underline X_1}{\underline{F}}$ and $y = \frac{X_0}{\underline{F}}$ the first equation reads
$$
x- \ln x = y - \ln y,
$$
and since $y = \frac{X_0}{\underline{F}} > 1$ an application of Theorem \ref{le:julelemma} gives
$$
x = z\left(y\right)  y e^{-y},
$$
where $z$ is a decreasing function of $y$.
As the same argument applies to the upper barrier we conclude
$$
z\left(\frac{X_0}{\underline{F}}\right){X_0} e^{-\frac{X_0}{\underline{F}}} = \underline X_1 < X_1 < \overline X_1 = z\left(\frac{X_0}{\overline{F}}\right){X_0} e^{-\frac{X_0}{\overline{F}}}.
$$
Moreover, since
$$
x = -W(-ye^{-y}),
$$
where $W$ is the principal branch of the Lambert $W$-function, we also have
$$
-\underline{F} W\left(-\frac{X_0}{\underline{F}}e^{-\frac{X_0}{\underline{F}}} \right) = \underline X_1 < X_1 < \overline X_1 = -\overline{F} W\left(\frac{X_0}{\overline{F}} e^{-\frac{X_0}{\overline{F}}}\right) .
$$
Finally, applying Theorem \ref{le:julelemma} for estimating $z$ yields the following result:

\begin{theorem}\label{eq:est-general}
Suppose that $T$ is a trajectory of system \eqref{eq:general} with initial condition $(X_0,S_0)$ where $X_0 > F(S)$ and $S_0$ satisfies $\varphi(S_0) < d/p$.
Suppose that, until $T$ intersects $S = S_0$ the next time, $T$ stays in a part of the state space where  \eqref{eq:ass_general} is satisfied.
Then, for the next intersection of trajectory $T$ with $S = S_0$ at $(X_1,S_0)$ it holds that
$$
-\underline{F} W\left(-\frac{X_0}{\underline{F}}e^{-\frac{X_0}{\underline{F}}} \right)  < X_1 <  -\overline{F} W\left(\frac{X_0}{\overline{F}} e^{-\frac{X_0}{\overline{F}}}\right),
$$
and
\begin{align*}
e^{-\frac{X_0}{\underline{F}}} < z_1  e^{-\frac{X_0}{\underline{F}}} < \frac{X_1}{X_0}  < z_2 e^{-\frac{X_0}{\overline{F}}}  < z_0 e^{\frac{X_0}{\overline{F}}} <  e^{1 + \frac{X_0}{\overline{F}}},
\end{align*}
where $z_1 = z_1\left(\frac{X_0}{\underline{F}}\right)$, $z_2 = z_2\left(\frac{X_0}{\overline{F}}\right)$ and $z_0 = z_0\left(\frac{X_0}{\overline{F}}\right)$.
\end{theorem}

\noindent
As a remark, we note that the accuracy of the estimates in Theorem \ref{eq:est-general} improves when $F(S) = H(S)/\varphi(S)$ obeys less variation,
i.e., when one can take a tighter interval in assumption \eqref{eq:ass_general}.

As an example, for which we easily find the above imposed assumptions satisfied
and thereby conclude Theorem \ref{eq:est-general},
we consider the system
\begin{align}\label{eqmainRM_n}
\frac{\mathrm ds}{\mathrm d\tau}&=(h(s)-x)s \notag\\
\frac{\mathrm dx}{\mathrm d\tau}&=m\, (s-\lambda )\,x
\end{align}
in $x,s\geq 0$, supposing parameters $m > 0$, $\lambda \in (0, 1)$ and where $h$ is given by
\begin{equation}\label{standardh}
h(s)=(1-s)(s+a).
\end{equation}
Any standard Rosenzweig-MacArthur predator-prey system can be transformed into a system of type \eqref{eqmainRM_n}.
In particular,  it is equivalent with system \eqref{eq:general} when \eqref{eq:RM} holds with $n = 1$,
which can be seen by introducing the non-dimensional quantities
\begin{align}\label{def-parameters}
\tau \,=\, \int \frac{r K}{A + S\left(t\right)} dt, \quad s &\,=\, \frac{S}{K}, \quad x \,=\, \frac{q X}{r K}, \quad a \,=\, \frac{A}{K},  \notag\\
m \,=\, \frac{p - d}{r} \quad &\text{and} \quad \lambda \,=\, \frac{d A}{(p-d) K}.
\end{align}

Comparing to the general system \eqref{eq:general} we identify
$S = s, X = x, F(S) = h(s), \varphi(S) = s,  p = m$, and $d = m \lambda$.
Hence, since $a < h(s)=(1-s)(s+a) < h(\lambda)$ when $s < \lambda$,
we obtain \eqref{eq:ass_general} satisfied in the region $s < \lambda$ with $\underline{F} = a$ and $\overline{F} = h(\lambda)$.
Observing that $s = \lambda$ is the isocline $x' = 0$ and $x = h(s)$ is the isocline $s' = 0$ we conclude, from Theorem \ref{eq:est-general},
the following estimate for the minimal predator biomass, $x_{\text{min}}$:

\begin{corollary}\label{eq:vanliga}
Consider a trajectory of system \eqref{eqmainRM_n} starting at $(x_{max},\lambda)$ with $x_{max} > h(\lambda)$.
Suppose that $m > 0$, $\lambda \in (0, 1)$ and that \eqref{standardh} holds.
Then the minimal predator biomass, $x_{\text{min}}$, satisfies
\begin{align*}
e^{-\frac{x_{\text{max}}}{a}} < z_1  e^{-\frac{x_{\text{max}}}{a}} < \frac{x_{\text{min}}}{x_{\text{max}}}  < z_2 e^{-\frac{x_{\text{max}}}{h(\lambda)}}  < z_0 e^{\frac{x_{\text{max}}}{h(\lambda)}} <  e^{1 + \frac{x_{\text{max}}}{h(\lambda)}},
\end{align*}
where $z_1 = z_1\left(\frac{x_{\text{max}}}{a}\right)$, $z_2 = z_2\left(\frac{x_{\text{max}}}{h(\lambda)}\right)$ and $z_0 = z_0\left(\frac{x_{\text{max}}}{h(\lambda)}\right)$.
\end{corollary}

We remark that by shrinking the region $s < \lambda$ into $s < \lambda^*$ for some $0 < \lambda^* < \lambda$ we have the better estimates
$a < h(s) < h(\lambda^*)$ and may, for a trajectory starting at $(x_0, \lambda^*)$, estimate its next intersection with the line $s = \lambda^*$ at point $(x_1, \lambda^*)$. In particular, we then get
Corollary \ref{eq:vanliga}
but with $x_{\text{max}} = x_0$, $x_{\text{min}} = x_1$ and $h(\lambda) = h(\lambda^*)$.
Furthermore, we observe that the upper estimate is good for small $\lambda^*$ but becomes less efficient when the value of $\lambda$ increases beyond $a$.
The bounds in
Corollary \ref{eq:vanliga}
will be used by the authors in \cite{LS22} for estimating the size of the unique limit cycle to system \eqref{eqmainRM_n} which is the global attractor when $2\lambda + a < 1$.


\subsection*{Coexistence of predators}

As a final remark we consider the following system similar to \eqref{eq:general} but allowing for $n \geq 1$ predators:
\begin{align}\label{eq:general_n}
\frac{\mathrm dS}{\mathrm dt} &= H(S) - \sum_{i=1}^{n} q_i \varphi_i(S) X_i, \notag\\
\frac{\mathrm dX_i}{\mathrm dt} &= p_i \varphi_i(S) X_i - d_i X_i, \quad i = 1, \dots, n,
\end{align}
where the non-negative variable $S$ represents the prey biomass,
the non-negative variables $X_i$ represent the predator biomasses,
$\varphi_i$ is non-decreasing,
$\varphi(0) = H(0) = 0$,
parameters $p_i$, $q_i$, $d_i$ are positive.

Following \cite[page 2]{sp} we assume $p_i > d_i$. If not, the corresponding predator will die.
Using the time change $\tau = rt$,
where $\tau$ is the new time, and the variable
changes $s = \frac{S}{K}, x_i = \frac{q_i}{r K} X_i$,
we transform,  when
\begin{align*}
H(S) = r S \left(1 - \frac{S}{K} \right) \quad \text{and}\quad \varphi_i(S) = \frac{S}{S + A_i},
\end{align*}
system \eqref{eq:general_n} to the system
\begin{align*}
\frac{\mathrm ds}{\mathrm d\tau} &= \left(1 - s  - \sum_{i=1}^{n} \frac{x_i}{s + a_i} \right) s, \notag\\
\frac{\mathrm dx_i}{\mathrm d\tau} &= m_i \frac{s-\lambda_i}{s + a_i}x_i, \quad i = 1, \dots, n,
\end{align*}
where
$$
a_i = \frac{A_i}{K}, \quad m_i = \frac{p_i - d_i}{r}, \quad \text{and}
\quad \lambda_i = \frac{d_i A_i}{K(p_i - d_i)}.
$$
%
Systems of this type have been studied e.g. in \cite{osipovAlushta,osipovIJBC,osipoveuler,sp}.
%
%
In particular, extinction and coexistence results for predators can be found in \cite{osipovAlushta,sp} from which we recall the following statement, 
giving sufficient condition for extinction: 
%
%

\begin{statement}
\cite[Statement 2]{sp} Let $L=\frac{ \lambda_i
(1-\lambda_j) }{ \lambda_j (1-\lambda_i) }$ and $\lambda_i>\lambda_j$. If
$a_j>\frac { a_i } { L+a_i (L-1) }$,
then the predator $i$ goes extinct.
\end{statement}

Anyhow the condition is quite far from necessary and it is possible to use
the results in this work to find sufficient conditions for the opposite, that is
for coexistence of predators. The proof of the statement and similar known proofs for extinction essentially only use the equations for $x_i$ and the properties of the functions $\varphi_i$, the equation for $s$
is only used to make the obvious restriction $s<1$.
If we consider the case of two predators we notice that in the two dimensional coordinate planes, where one predator is abscent, there can be cycles like in the standard Rosenzweig-MacArthur system.
The instability of one or two of these cycles can be used to get sufficient conditions for coexistence, and thus we conjecture that nice estimates we have produced in this work could be used for getting important coexistence conditions contradicting the known exclusion principle.

Another interesting problem is arising in connection to the bifurcation for this system examined in \cite{nico}. Here they under certain conditions conclude cyclic coexistence of the predators for parameters near to the case $\lambda_1=\lambda_2$ where there is a cycle only in one of the coordinate planes. We conjecture that our estimates for this cycle can be used to prove its instability and thus coexistence for a parameter range after the bifurcation.

We observe that in all these cases we suppose to use the results only far beyond the
bifurcation to cycle of the equilibrium when the cycles are big and parameters $a_i$ and
$\lambda_i$ are small.

It occurs that the behaviour for small prey biomass on the cycle can play an important role in determining the stability,
see \cite{osipoveuler} and references therein,
and thus estimates of the limit cycle for small prey, such as those in Theorem \eqref{eq:est-general}, obtained by Theorem \ref{le:julelemma}, may be useful. \\


\noindent
{\bf Acknowledgment}\\
The authors would like to thank two anonymous reviewer for valuable comments and suggestions which
really helped us to improve this work.



\begin{thebibliography}{99}

\bibitem{ACL24}
{\AA}hag P., Czyz R.,  Lundow P-H.
{\it On a generalised Lambert W branch transition function arising from p, q-binomial coefficients.}
Applied Mathematics and Computation 462, (2024): 128347.


\bibitem{nico}
Alejandro L\'{o}pez-Nieto, Phillipo Lappicy, Nicola Vassena,
Hannes Stuke, Jia-Yuan Dai.
{\it Hybrid Bifurcations and Stable Periodic
Coexistence for Competing Predators.}
https://arxiv.org/abs/2310.19604

\bibitem{hej11}
Aravindakshan A., Ratchford B.
{\it Solving share equations in logit models using the lambertW function}.
Review of Marketing Science, 9(1): (2011), 1-17.

\bibitem{BPLPCS00}
Barry D. A., Parlange J. Y., Li L., Prommer H., Cunningham  C. J.,  Stagnitti F.
{\it Analytical approximations for real values of the Lambert W-function.}
Mathematics and Computers in Simulation, 53(1-2), (2000): 95-103.

%
\bibitem{C81}
Cheng K.-S.
{\it Uniqueness of a limit cycle for a predator-prey system}.
SIAM Journal on Mathematical Analysis 12.4 (1981): 541-548.
%
\bibitem{CV98}
Clanet C., Villermaux E.
{\it An analytical estimate of the period for the delayed logistic application and the Lotka-Volterra system.}
The European Physical Journal B-Condensed Matter and Complex Systems, 6 (1998): 529-536.
%
%
%
%
\bibitem{GS07}
Grozdanovski T., Shepherd J.
{\it Approximating the periodic solutions of the Lotka-Volterra system.}
ANZIAM Journal, 49 (2007): C243-C257.
%
%
%
\bibitem{HH08}
Hoorfar A.,  Hassani M.
{\it Inequalities on the Lambert W function and hyperpower function.}
J. Inequal. Pure and Appl. Math 9.2, (2008): 5-9.


%
%
\bibitem{HS09}
Hsu, S.-B., Shi, J.:
{\it Relaxation oscillator profile of limit cycle in predator-prey system}.
Discrete and continuous dynamical systems B 11, no. 4 (2009): 893-911.
%
%
%
\bibitem{IDM20}
Ito H. C., Dieckmann U., Metz J. A.
{\it Lotka-Volterra approximations for evolutionary trait-substitution processes.}
Journal of Mathematical Biology, 80(7), (2020): 2141-2226.
%
%
%
\bibitem{L89}
Lindstr\"om T.
{\it A generalized uniqueness theorem for limit cycles in a predator-prey system.}
Acta Academiae Aboensis, Series B: Mathematica et Physica, 49(2), (1989).
%
\bibitem{L94}
Lindstr\"om T.
{\it Why do rodent populations fluctuate?: stability and bifurcation analysis of some discrete and continuous predator-prey models}
Doctoral dissertation, Luleå tekniska universitet (1994).
%
\bibitem{LS18}
Lundstr\"om N. L. P., S\"oderbacka G. J.,
{\it Estimates of size of cycle in a predator-prey system},
{Differential Equations and Dynamical Systems}, (2018): 1-29. 
%
\bibitem{LS22}
Lundstr\"om N. L. P., S\"oderbacka G. J.,
{\it Estimates of size of cycle in a predator-prey system II},
https://arxiv.org/abs/2311.16735
%
%
%
%
\bibitem{M05}
Morales D. A.,
{\it Exact expressions for the range and the optimal angle of a projectile with linear drag}.
Canadian Journal of Physics {\bf 83}, (2005): 67-83.

%
%
\bibitem{MR87}
Murty K. N.,  Rao D. V. G.
{\it Approximate analytical solutions of general Lotka-Volterra equations.}
J. Math. Anal. Applic., 122(2), (1987): 582-588.
%
\bibitem{osipovAlushta}
Osipov A. V.,  S{\"o}derbacka G. J.,
{\it Extinction and coexistence of predators}.
{Dinamicheskie Sistemy}, {\bf 6}, 1, (2016): 55-64.
http://dynsys.cfuv.ru/
%
\bibitem{osipovIJBC}
Osipov A. V.,   S{\"o}derbacka G. J.,
{\it Poincar\'{e}  map construction for some classic two predators-one prey systems}.
{Internat. J. Bifur. Chaos Appl. Sci. Engrg},{\bf 27}, (8), (2017): 1750116, 9 pp.
%
\bibitem{osipoveuler}
Osipov A. V.,   S{\"o}derbacka G. J.,
{\it Review of results on a system of type many predators-one prey}.
{Nonlinear Systems and Their Remarkable Mathematical Structures} (2018): 520-540. CRC Press.
%
\bibitem{PY04}
Packel E. W., D. S. Yuen,
{\it Projectile Motion with Resistance and the Lambert W Function}.
The College Mathematics Journal {\bf 35}, (2004): 337.

%
%
%
%


\bibitem{hej14}
Sharma S., Shokeen P., Saini B., Sharma S., Chetna, Kashyap J., Guliani R., Sharma S., Udaibir, Khanna M., Jain A., Kapoor A.
{\it Exact analytical solutions of the parameters of different generation real solar cells using Lambert W-function: A Review Article.}
Invertis Journal of Renewable Energy, 4(4), (2014): 155-194.

\bibitem{sp}
S\"{o}derbacka G. J., Petrov A. S.
Review on the behaviour of a many predator - one prey system.
 {\it Dinamicheskie sistemy},  9(37), No3, (2019), 273-288.

\bibitem{S09}
Stewart S. M.
{\it On certain inequalities involving the Lambert W function}.
Journal of Inequalities in Pure and Applied Mathematics {\bf 10} (2009).

\bibitem{S06}
Stewart S. M.
{\it An analytic approach to projectile motion in a linear resisting medium}.
International Journal of Mathematical Education in Science and Technology {\bf 37},  (2006), 411-431.


%
%
%
\bibitem{WW04}
Warburton R. D. H., Wang J.,
{\it Analysis of asymptotic projectile motion with air resistance using the Lambert W function.}
American Journal of Physics {\bf 72}, (2004): 1404-1407.

\end{thebibliography}
\end{document}